\documentclass[reqno,11pt]{amsart}

\usepackage{amsmath,amssymb,amsthm,mathrsfs}
\usepackage{epsfig,color}
\usepackage[latin1]{inputenc}

\usepackage[numbers]{natbib}

 \allowdisplaybreaks
\numberwithin{equation}{section}

\def\R{\mathbb R}
\def\N{\mathbb N}

\def\e{\varepsilon}
\def\s{\sigma}

\def\Div{{\rm div}\,}

\def\l{\lambda}
\def\g{\gamma}
\def\k{\kappa}

\def\de{\delta}

\def\pa{\partial}

\renewcommand{\a}{\alpha}

\renewcommand{\l}{\lambda}
\renewcommand{\L}{\Lambda}

\newcommand{\n}{\nabla}

\renewcommand{\Div}{{\rm div \,}}

\theoremstyle{plain}
\newtheorem{theorem}{Theorem}[section]

\newtheorem{remark}[theorem]{Remark}

\title[Metrics with almost constant positive scalar curvature]{A quantitative analysis of metrics on $\R^n$ with almost constant positive scalar curvature, with applications to fast diffusion flows}

\author{G. Ciraolo}
\address{ Dipartimento di Matematica e Informatica,
Università di Palermo, Via Archirafi 34, 90123 Palermo, Italy}
\email{giulio.ciraolo@unipa.it}

\author{A. Figalli}
\address{Mathematics Department, ETH Z\"urich,
R\"amistrasse 101, CH-8092 Z\"urich, Switzerland}
\email{alessio.figalli@math.ethz.ch}

\author{F. Maggi}
\address{Abdus Salam International Center for Theoretical Physics,
Strada Costiera 11, I-34151, Trieste, Italy.
On leave from the University of Texas at Austin}
\email{fmaggi@ictp.it}

\begin{document}

\begin{abstract}
{\rm We prove a quantitative structure theorem for metrics on $\R^n$ that are conformal to the flat metric, have almost constant positive scalar curvature,
and cannot concentrate more than one bubble. As an application of our result, we show a quantitative rate of convergence in relative entropy for a fast diffusion equation in $\R^n$ related to the Yamabe flow.}
\end{abstract}

\maketitle

\section{Introduction}

\subsection{The prescribed scalar curvature problem on the sphere}
Given a $n$-dimensional Riemannian manifold $(M,g_0)$, $n\ge 3$, the problem of finding a metric $g$ conformal to $g_0$ whose scalar curvature $R_g$ is equal to a prescribed function $R$ boils down to showing the existence of a positive solution $u$ to the nonlinear PDE
\begin{equation}
  \label{prescribed scalar curvature on M}
  -\Delta_{g_0}u+\frac{n-2}{4(n-1)}\,R_{g_0}\,u=\frac{n-2}{4(n-1)}\,R\,u^p
\end{equation}
where $\Delta_{g_0}=\Div(\nabla_{g_0}\cdot )$ and $R_{g_0}$ denote respectively the Laplace-Beltrami operator
and the scalar curvature of $(M,g_0)$, and
\[
p=\frac{n+2}{n-2}=2^\star-1\,,\qquad 2^\star=\frac{2n}{n-2}\,.
\]
Indeed, if $u$ solves \eqref{prescribed scalar curvature on M}, then the metric $g=u^{p-1}\,g_0$ satisfies $R_g=R$.

When $(M,g_0)$ is the round sphere then $R_{g_0}=n(n-1)$ and \eqref{prescribed scalar curvature on M}
can be read on $\R^n$ by means of the stereographic projection.
More precisely,
consider the inverse stereographic projection $F:\R^n\to \mathbb{S}^n$
defined by
\begin{equation}
\label{eq:F}
F(x)=\Bigl(\frac{2x}{1+|x|^2},\frac{|x|^2-1}{1+|x|^2}\Bigr)\,.
\end{equation}
Then $v:\mathbb S^n\to \R$ solves \eqref{prescribed scalar curvature on M} if and only if
$$
u(x)=\Bigl(\frac{2}{1+|x|^2}\Bigr)^{(n-2)/2}v(F(x))
$$
solves
\begin{equation}
  \label{prescribed scalar curvature equation in Rn}
  -\Delta u=K\,u^p\quad\mbox{on $\R^n$}\,,
\end{equation}
where $K(x)=\frac{n-2}{4(n-1)}\,R(F(x))$.

When looking for solutions of \eqref{prescribed scalar curvature equation in Rn},
it is natural to impose that $u$ satisfies
\begin{equation}
  \label{positive solution with finite dirichlet}
  \mbox{$u>0$ on $\R^n$ }\qquad \text{and}\qquad \int_{\R^n}|\nabla u|^2<\infty\,.
\end{equation}
In general, for a given function $K$ (and even when $K$ is just a small perturbation of a constant) there may exist no solution to \eqref{prescribed scalar curvature equation in Rn}--\eqref{positive solution with finite dirichlet} (see \cite{kazdanwerner74}), and indeed there is a vast literature dedicated to finding necessary and sufficient conditions on $K$ in order to guarantee the solvability of \eqref{prescribed scalar curvature equation in Rn}--\eqref{positive solution with finite dirichlet}, see for example \cite{moser71,ni82,dingni85,escobarschoen86,chen89,changyang91,bahricoron91,bianchiegnell93,YYLiPARTI,YYLiPARTII,bahriduke96,ChenLi97,ambrosettiazoreroperal99,ambromalchiodiJDE01,LinLiu09,LeungZhou15}. At the same time, when $K$ is constantly equal to some $\k>0$, the problem is completely rigid. Indeed, by \cite{obata,gidasninirenberg}, if $u$ solves \eqref{prescribed scalar curvature equation in Rn}--\eqref{positive solution with finite dirichlet} with $K\equiv\k>0$, then there exist $\l>0$ and $z\in\R^n$ such that
\[
u(x)=\l^{(n-2)/2}\,v_\k(\l(x-z))\qquad\forall\, x\in\R^n\,,
\]
where
\begin{equation}
  \label{standard metric of the sphere}
  v_\k(x)=\Big(\frac{n(n-2)}{\k}\Big)^{(n-2)/4}\,\frac{1}{(1+|x|^2)^{(n-2)/2}}\qquad\forall\, x\in\R^n\,.
\end{equation}

\subsection{Main result}
The goal of this paper is to give a quantitative description of solutions $u$ to the {\it prescribed scalar curvature equation} \eqref{prescribed scalar curvature equation in Rn}--\eqref{positive solution with finite dirichlet} in the regime when $K$ is close (in a suitable sense) to a positive constant.

In order to identify the natural space in which the distance of $K$
to a constant should be measured, we make the following observation:
$$
-\Delta u =K\,u^p \qquad \Longleftrightarrow\qquad
\int_{\R^n} \n u \cdot \n\varphi=\int_{\R^n} K\,u^p\varphi \qquad \forall\,\varphi \in \dot W^{1,2}(\R^n),
$$
where $\dot W^{1,2}(\R^n)$ denotes the closure of $C^\infty_c(\R^n)$ with respect to the norm $\|\nabla\cdot\|_{L^2(\R^n)}$.
Since
$$
\varphi \in \dot W^{1,2}(\R^n)\hookrightarrow L^{2^*}(\R^n),
$$
$\int_{\R^n} K\,u^p\varphi$ is well-defined provided $K\,u^p$  belongs to the dual of $L^{2^*}(\R^n)$, namely $L^{2n/(n+2)}(\R^n).$

This suggests the following definition:
\begin{equation}
  \label{deficit}
  \de(u):=\|K\,u^p-K_0(u)\,u^p\|_{L^{2n/(n+2)}}=\Big(\int_{\R^n}\Big|\frac{\Delta u}{u^p}+K_0(u)\Big|^{2n/(n+2)}\,u^{2^\star}\Big)^{(n+2)/2n}\,,\end{equation}
where
\begin{equation}
  \label{eq:K0}
K_0(u):=\frac{\int_{\R^n} K u^{2^\star}}{\int_{\R^n}u^{2^\star}}= \frac{\int_{\R^n} \frac{-\Delta u}{u^p}u^{2^\star}}{\int_{\R^n}u^{2^\star}}=
\frac{\int_{\R^n}|\nabla u|^2}{\int_{\R^n}u^{2^\star}}.
\end{equation}
Hence, the question becomes: if $\de(u)$ is small, can we say that $u$
is close to a translation/dilation of $v_{K_0(u)}$?

A negative answer is given by the following simple example:
given a function $v:\R^n\to\R$, set
\[
v[z,\l](x):=\l^{(n-2)/2}\,v(\l(x-z))\,,\qquad x\in\R^n\,,
\]
and consider
$$
 u=\sum_{i=1}^m v_1[z_i,\l_i],
 $$
 where the functions $v_1[z_i,\l_i]$ are supported far away from each other.
Then
 $$
 -\Delta u=-\sum_{i=1}^m \Delta v_1[z_i,\l_i]
 =\sum_{i=1}^m v_1[z_i,\l_i]^p=K \biggl(\sum_{i=1}^m v_1[z_i,\l_i]\biggr)^p=Ku^p,
 $$
 with
 $$
 K=\frac{\sum_{i=1}^m v_1[z_i,\l_i]^p}{\left(\sum_{i=1}^m v_1[z_i,\l_i]\right)^p},
 $$
 and it is easy to check that by taking the point $z_i$ sufficiently far from each other one can make $\|K\,u^p-\,u^p\|_{L^{2n/(n+2)}}$ arbitrarily small.

As shown by Struwe's \cite{struwe} (see also \cite[Theorem 3.3]{hebey}),
this bubbling phenomenon is the only possible ``bad'' case. More precisely, whenever $\delta(u)$ is small, $u$ is close to a sum of bubbles as above.
Luckily,
in many applications, this phenomenon can be avoided by some preliminary study
of the PDE under investigation, and one can usually localize the problem in a suitable way so that, in the region under investigation, the ``energy'' of $u$ is strictly less than the energy of two bubbles.
Hence, we shall focus on the latter situation.

Before stating our result, we recall the definition of the Sobolev constant on $\R^n$,
\begin{equation}
  \label{sobolev constant}
  S=\inf\Big\{\frac{\|\nabla v\|_{L^2(\R^n)}}{\|v\|_{L^{2^\star}(\R^n)}}:v\ne 0\,,\ |\{|v|>t\}|<\infty\ \ \forall\, t>0\Big\}\,.
\end{equation}
By \cite{aubin,talenti,cenv}, the family of functions $\{v_\k[z,\l]\}_{\k\,,\l\,,z}$ corresponds to the minimizers in \eqref{sobolev constant}.
Hence, since $-\Delta v_\k=\k\,v_\k^p$, one can easily check that
\begin{equation}
  \label{talenti opt}
  \int_{\R^n}|\nabla v_\k|^2=S^2\,\Big(\int_{\R^n}v_\k^{2^\star}\Big)^{2/2^\star}=\frac{S^n}{\k^{(n-2)/2}}\qquad\int_{\R^n}v_\k^{2^\star}=\frac{S^n}{\k^{n/2}}\,.
\end{equation}

To simplify the notation (and because for most applications this is not a real restriction), we shall assume that the energy of $u$ is bounded by $3/2$ the energy of a single bubble. Of course $3/2$ does not play any essential role, and the proof holds when $3/2$ is replaced by any constant strictly less than $2$.
Also, the assumption $K_0(u)=1$ is not restrictive, since it can always be guaranteed by rescaling $u$.

\begin{theorem}
  \label{thm main}
  Given $n\ge 3$, there exists $C_0=C_0(n)>0$ with the following property. Let $u\in C^\infty(\R^n)\cap \mathring{H}^1(\R^n)$
  be a positive function satisfying
  \[
  K_0(u)=1\qquad \text{and}\qquad \int_{\R^n}|\nabla u|^2\le \frac{3}{2}\,S^n\,.
  \]
  Then there exist $z\in \R^n$ and $\l \in (0,\infty)$ such that
  \[
  u=\,v_1[z,\l]+\rho\,,
  \]
  where
  \begin{equation}\label{tesi thm main 1}
  \|\nabla\rho\|_{L^2(\R^n)}
  \le C_0\de(u).
  \end{equation}
\end{theorem}

\begin{remark}
  {\rm Theorem \ref{thm main} is easily seen to be optimal. Indeed, set $v:=v_1[0,1]$, let $\e>0$ be small, and consider $u_\e:=v_1+ \e \phi$ for some $\phi\in C^\infty_c(\R^n),$ with $\phi \geq 0$.
Then
$$
\Delta u_\e = \Delta v +\e \Delta \phi=v^p+\e \Delta \phi=K_\e u^p
$$
with
$$
K_\e:=\frac{v^p+\e \Delta \phi}{u_\e^p}=1+\frac{[v^p-(v+\e \phi)^p]+\e \Delta \phi}{u_\e^p}=1+O(\e u_\e^{-p}).
$$
This shows that
\begin{equation}
\label{eq:Ke}
\|K_\e u_\e^p - u_\e^p\|_{L^{2n/(n+2)}}=O(\e).
\end{equation}
Also, since 
$$
K_0(u_\e)=\frac{\int_{\R^n}K_\e u_\e^{2^\star}}{\int_{\R^n}u_\e^{2^\star}},
$$
by \eqref{eq:Ke} and H\"older inequality,
$$
|K_0(u_\e)-1|\leq \frac{|\int_{\R^n}(K_\e-1) u_\e^{2^\star}|}{\int_{\R^n}u_\e^{2^\star}} \leq 
\biggl(\frac{\int_{\R^n}|K_\e-1|^{2n/(n+2)} u_\e^{2^\star}}{\int_{\R^n}u_\e^{2^\star}} \biggr)^{(n+2)/2n}=O(\e).
$$
Combining this bound with \eqref{eq:Ke} we get $\delta(u_\e)\leq C\e \simeq \|\nabla (\e\phi)\|_{L^2(\R^n)}$, which proves the optimality of our result. 
  }

\end{remark}

As an application of Theorem \ref{thm main},
we investigate the behavior of solutions to
a fast diffusion equation in $\R^n$ related to the Yamabe flow.

\subsection{Convergence to equilibrium: a fast diffusion
equation related to the Yamabe flow}
\label{sect:intro FD}
Given $m \in (0,1)$, the Cauchy problem for the fast diffusion equation
is written as
\begin{equation}
\label{eq:FDm}
\frac{d}{dt} u =\Delta(u^m)\qquad \text{in }(0,\infty)\times \R^n\,.
\end{equation}
Assuming that the initial datum $u_0$ is nonnegative and fastly decaying at infinity,
it is well-known that solutions to \eqref{eq:FDm} are smooth and positive for all
times if $m>m_c=(n-2)/n$, while they vanish in finite time if $m \leq m_c$.
There is a huge literature on the subject, and we refer the interested reader to the monograph
 \cite{vazquezBook} for a comprehensive overview and more references.

A case of a special interest corresponds to the choice $m=(n-2)/(n+2)$,
where the equation for $u$ is equivalent to the (not volume-preserving) Yamabe flow
$$
\frac{d}{dt}g=-R_g\,g
$$
 for the metric $g_{ij}(t)=u(t)^{p-1}\,dx_i\,dx_j$
 (recall that $p=2^\star-1$).
Hence, we consider the Cauchy problem
\begin{equation}
\label{eq:FD}
\frac{d}{dt} u =\Delta(u^m)\qquad \text{in }(0,\infty)\times \R^n\,,\qquad m=\frac{n-2}{n+2}=\frac1p\,,
\end{equation}
with a continuous initial datum $u_0\geq 0$ satisfying $u_0(x)=O(|x|^{-(n+2)})$ as $|x|\to \infty$.
Because of its geometric relevance, this equation has received a lot of attention.
In particular, as shown in \cite[Theorem 1.1]{DPS} (see also \cite{chow,yerugang}), under the above assumptions on $u_0$ there exists a vanishing time $T=T(u_0)>0$, a point $z \in \R^n$, and a number $\l>0$, such that $u\equiv 0$ for $t \geq T$ and
\begin{equation}
\label{eq:convergence}
\biggl\|\frac{(T-t)^{-1/(1-m)}u(t,\cdot)}{v_{1/(1-m)}[z,\l]^{1/m}} -1\biggr\|_{L^\infty(\R^n)} \to 0\qquad \text{as }t \to T^-.
\end{equation}
Later on, in \cite{DaskSes,BBDGV,BDGV,BGV10}, the authors
investigated the asymptotic of solutions for all values of $m \in (0,1)$,
and proved both qualitative and quantitative convergence results under the assumption that
the initial datum is trapped in between two Barenblatt solutions
with the same extinction time.

As observed in \cite{DaskSes} (see also \cite{DKS13} for a different but related analysis), this trapping assumption on
the initial datum is very restrictive in our setting,
as it completely misses
the picture given by \cite{DPS} and gives rise to an extinction profile different from
the one in \eqref{eq:convergence}, which is believed to be the correct one for
``most'' initial data.
Indeed, as shown in \cite[Theorem 1.4]{DaskSes},
there exists a large class of initial data of the form
$$
u_0(x)=\frac{C_0}{|x|^{(n+2)/2}}(1+o(1))\qquad \text{as $|x|\to \infty$}\,,\quad C_0>0\,,
$$
that do not satisfy the trapping assumption and
whose solutions behave as follows:
there exists a time $t_0\in (0,T)$ (which can be explicitly computed in terms of $C_0$
and is given by $t_0=\frac{C_0(n+2)}{(n-2)^2}$)
 such that
the behavior of the solution is governed by
a Barenblatt profile $B(t,\cdot)$ up to the time $t_0$ when the Barenblatt vanishes. Then the solution develops
a singularity at $t=t_0$, and it satisfies $u(t,x)=O(|x|^{-(n+2)})$ as $|x|\to \infty$
for all $t>t_0$
(see \cite[Theorem 1.4(ii)]{DaskSes}).
In particular, as observed in \cite[Corollary 1.5]{DaskSes}, this allows one to apply \cite[Theorem 1.1]{DPS}
and deduce that $u(t,x)$ exhibits the vanishing profile of a sphere,
as shown by \eqref{eq:convergence}.

By exploiting our Theorem \ref{thm main}, we can improve the convergence result in
\eqref{eq:convergence}
and obtain a quantitative rate of convergence under the same assumptions as in \cite{DPS}. More precisely,
in section \ref{sect:FD} we prove the following:

\begin{theorem}
\label{thm:FD}
Let $u$ be a solution to the fast diffusion equation \eqref{eq:FD}
starting from a non-negative continuous initial datum $u_0$ satisfying
 $u_0(x)=O(|x|^{-(n+2)})$ as $|x|\to \infty$.
 Let $T=T(u_0)>0$ denote the vanishing time of $u$.
Then there exist $z \in \R^n$, $\l>0$, and a dimensional constant $\k(n)>0$ such that
$$
\biggl\|\frac{(T-t)^{-1/(1-m)}u(t,\cdot)}{v_{1/(1-m)}[z,\l]^{1/m}} -1\biggr\|_{L^\infty(\R^n)} \leq C_*\,(T-t)^{\k(n)}\qquad \forall\, 0<t<T\,,
$$
where $C_*>0$ is a constant depending on the initial datum $u_0$.
\end{theorem}

\bigskip

\noindent {\bf Acknowledgment:} The authors wish to thank Manuel Del Pino for suggesting the possibility of improving \cite[Theorem 1.1]{DPS}. This work has been done while GC was visiting the University of Texas at Austin under the support of NSF FRG Grant DMS-1361122,
the Oden Fellowship at ICES, the GNAMPA of the Istituto Nazionale di Alta Matematica (INdAM), and the FIR project 2013 ``Geometrical and Qualitative aspects of PDE''.
 AF was partially supported by NSF Grants DMS-1262411 and DMS-1361122, and FM was supported by NSF Grants DMS-1265910
 and DMS-1361122.

\section{Proof of Theorem \ref{thm main}}\label{proof of the main theorem} This section is devoted to the proof of Theorem \ref{thm main}. Recall the notation
\[
n\ge 3\,,\qquad 2^\star=\frac{2n}{n-2}\,,\qquad p=2^\star-1=\frac{n+2}{n-2}=1+s\,,\qquad s=p-1=\frac{4}{n-2}\,.
\]
We begin by observing that it is enough to prove the theorem for $u:\R^n\to \R$ satisfying
\[
u\in C^\infty(\R^n)\cap \mathring{H}^1(\R^n)\,,\qquad \mbox{$u>0$ on $\R^n$}\,,
\]
and such that
\begin{equation}
  \label{assumption 1}
  K_0(u)=1\,,\qquad \int_{\R^n}|\nabla u|^2\le \frac{3}{2}\,S^n\,,\qquad \de=\de(u)\le\de_0\,,
\end{equation}
for a suitably small constant $\de_0=\de_0(n)$.

Indeed we note that if $\delta(u)>\delta_0$, then the theorem is trivially true simply by choosing
$\lambda=1$, $z=0$,
setting
$$
\rho=u-v_1\,,
$$
and then simply choosing $C_0=C_0(n)$ large enough.

\medskip

\noindent {\it Step one}.
Thanks to \cite{struwe}, for a constant $\e_0=\e_0(n)>0$ to be fixed later on, we can choose $\de_0$ depending on $\e_0$ in such a way that there exist $z\in \R^n$, $\l,\a>0$, such that
\begin{equation}  \label{struwe conclusion proof}
\begin{split}
\Big\|\nabla u- \alpha \nabla v_1[z,\l]\Big\|_{L^2(\R^n)}\le\e_0\,,\\
|\alpha -1|\leq \e_0.
\end{split}
\end{equation}
Without loss of generality, the parameters  $z\in\R^n$ and $\l,\a\in(0,\infty)$ can be chosen in such a way that
\begin{equation}
  \label{minimal choice}
  \Big\|\nabla u-\alpha \nabla v_1[z,\l]\Big\|_{L^2(\R^n)}=\min_{\,w\in\R^n,\,\mu,a>0}\Big\|\nabla u- a \nabla v_1[w,\mu]\Big\|_{L^2(\R^n)}\,.
\end{equation}
In particular, if we set
\[
\begin{split}
  &\hspace{3cm}
   \rho=u-\sigma\,,\qquad \sigma=\a U\,,
  \\
  &U=v_1[z,\l]\,,\qquad V=\frac{\pa v_1[w,\mu]}{\pa\mu}\bigg|_{w=z\,,\mu=\l}\,,\qquad W^j=\frac{\pa v_1[w,\mu]}{\pa w_j}\bigg|_{w=z\,,\mu=\l}\,,
\end{split}
\]
then by \eqref{minimal choice} we find that
\begin{equation}
  \label{rho orthogonal to}
  \int_{\R^n}\nabla U\cdot\nabla\rho=\int_{\R^n}\nabla V\cdot\nabla\rho=\int_{\R^n}\nabla W^j\cdot\nabla\rho=0\,,\qquad \forall \,1\le j\le n\,.
\end{equation}
Also, the first bound in \eqref{struwe conclusion proof} gives
\begin{equation}
  \label{rho small in H1}
  \|\nabla\rho\|_{L^2(\R^n)}\le\e_0\,.
\end{equation}
By a spectral analysis argument (see, e.g., the appendix to \cite{bianchiegnell}), \eqref{rho orthogonal to} implies that
\begin{equation}
  \label{spectral gap lowerbound}
  \int_{\R^n}|\nabla\rho|^2\ge \Lambda\,\int_{\R^n}U^{p-1}\rho^2 
\end{equation}
where $\Lambda=\Lambda(n)$ is such that
$$
  \Lambda>p\,.
$$
\bigskip

  \noindent {\it Step two}. Having in mind to exploit $\Lambda>p$,
we now test the equation   $-\Delta u=K\,u^p$ with $\rho$, and using $\int_{\R^n}\nabla\rho\cdot\nabla U=0$ we get
\begin{equation}
    \label{single bubble proof 1 y}
      \int_{\R^n}|\nabla\rho|^2=\int_{\R^n}K\,u^p\,\rho=\int_{\R^n}u^p\,\rho+\int_{\R^n}(K-1)\,u^p\,\rho\,.
  \end{equation}
  Since $u=\sigma+\rho$, a Taylor expansion yields
   \begin{equation}
    \label{intermedia}
      \int_{\R^n}u^p\,\rho=\int_{\R^n}\s^p\rho+p\,\int_{\R^n}\s^{p-1}\rho^2
  +{\rm O}\Big(\int_{\R^n}|\nabla\rho|^2\Big)^{1+\g}\,,
  \end{equation}
  where
    \[
  \g=\min\Big\{\frac12,\frac{2}{n-2}\Big\}\,.
  \]
  From $\s^p=\a^p U^p=-\a^p \Delta U$ and using again that $\int_{\R^n}\nabla\rho\cdot\nabla U=0$ we get
  $$
  \int_{\R^n}|\nabla\rho|^2-p\,\int_{\R^n}\s^{p-1}\rho^2=\int_{\R^n}(K-1)\,u^p\,\rho+{\rm O}\Big(\int_{\R^n}|\nabla\rho|^2\Big)^{1+\g}.
  $$
  Note that,
  by H\"older inequality (recall that $p=2^\star/(2^{\star})'$ and $(2^\star)'=2n/(n+2)$)
  \begin{equation}
    \label{new}
     \Big|\int_{\R^n}(K-1)\,u^p\,\rho\Big|\le\de(u)\,\|\rho\|_{L^{2^\star}(\R^n)}
  \le S\,\|\nabla\rho\|_{L^2(\R^n)}\,\de(u)\,.
  \end{equation}
  Also, recalling \eqref{spectral gap lowerbound},
  $$
  \int_{\R^n}\s^{p-1}\rho^2=\a^{p-1}\int_{\R^n}U^{p-1}\rho^2\leq \frac{\a^{p-1}}{\L}\int_{\R^n}|\nabla\rho|^2.
  $$
  Hence
  $$
  \Bigl(1-\frac{\a^{p-1}p}{\L}\Bigr)\int_{\R^n}|\nabla\rho|^2\leq S\,\|\nabla\rho\|_{L^2(\R^n)}\,\de(u)+{\rm O}\Big(\int_{\R^n}|\nabla\rho|^2\Big)^{1+\g}.
  $$
  Since $\L>p$ and $|\a-1|\leq \e_0$, choosing $\de(u)$ small enough we can ensure that $1-\frac{\a^{p-1}p}{\L} \geq c_0>0$ for some dimensional constant $c_0$, and we get
  $$
  c_0\int_{\R^n}|\nabla\rho|^2 \leq S\,\|\nabla\rho\|_{L^2(\R^n)}\,\de(u)+{\rm O}\Big(\int_{\R^n}|\nabla\rho|^2\Big)^{1+\g}.
  $$
  Since $\int_{\R^n}|\nabla\rho|^2$ is smaller than $\e_0$, we can also reabsorb the last term to conclude that $\|\nabla\rho\|_2 \leq C\,\delta(u)$.

  \bigskip

  \noindent {\it Step three}. We now quantitatively control $|\alpha-1|$.
  To this aim, we observe that assumption $K_0(u)=1$ is equivalent to
  $$
  \int_{\R^n}|\n u|^2=\int_{\R^n}u^{2^\star}
  $$
  (see \eqref{eq:K0}).
Note that, since $\|\nabla\rho\|_2 \leq C\,\delta(u)$ (by Step 2),
$$
  \int_{\R^n}|\n u|^2=  \int_{\R^n}|\n \s|^2+  \int_{\R^n}|\n \rho|^2
  =\a^2S^n+O(\delta^2).
$$
On the other hand, recalling that $\int_{\R^n}\s^{p}\rho=0$
(cp. with Step two)
$$
\int_{\R^n}u^{2^\star}=\int_{\R^n}\s^{2^\star}
+p\int_{\R^n}\s^{p}\rho+O(\de^2)=\alpha^{2^*}S^n+O(\delta^2).
$$
Comparing these expressions, we immediately deduce that
$$
\alpha^{2^*}S^n=\a^2S^n+O(\delta^2)\,,
$$
hence $|\a-1|\leq O(\de^2)$.
Thus, if we set $\rho':=\rho+(\a-1)U$, we proved that
$$
u=U+\rho'
$$
where $\|\n\rho'\|_2\leq C_0\delta(u)$, as desired.\qed

\section{An application to fast diffusion equations: Proof of Theorem \ref{thm:FD}}
\label{sect:FD}

This section is devoted to the proof of Theorem \ref{thm:FD}. Thus, we consider a
continuous initial datum $u_0\ge0$ satisfying $u_0(x)=O(|x|^{-(n+2)})$ as $|x|\to \infty$,
and $u$ a solution of \eqref{eq:FD} with $u(0,\cdot)=u_0$.
As explained in Section \ref{sect:intro FD},
under these assumptions the qualitative convergence result \eqref{eq:convergence}
holds,
and our goal is to quantify the rate of convergence.

We first notice that \eqref{eq:convergence} can be restated as follows
(see \cite[Theorem 1.1]{DPS}): there exists a vanishing time $T>0$ (depending on $u_0$) such that $u\equiv 0$ for $t \geq T$ and
\begin{equation}
\label{eq:vanishing}
\frac{u(t,x)}{(T-t)^{1/(1-m)}}=v_{1/(1-m)}[z_\infty,\l_\infty]^{1/m}+\theta(t,x)\quad \text{for }t<T\,,
\end{equation}
where $\lambda_\infty>0$, $z_\infty \in \R^n$, and $\theta$ satisfy
\begin{equation}
  \label{delpinosaez 1}
  \sup_{x \in \R^n}(1+|x|^{n+2})|\theta (t,x)|\to 0\qquad \text{as }t \to T^-\,.
\end{equation}
For proving Theorem \ref{thm:FD} we need to show the existence of $\k(n)>0$ such that
\begin{equation}
\label{eq:to prove FD}
\sup_{x \in \R^n}(1+|x|^{n+2})|\theta (t,x)|\leq C(T-t)^{\k(n)}\qquad \text{for }0<t<T\,.
\end{equation}
Following \cite{DPS}, define
\[
w(s,x)=\frac{u(t,x)^m}{(T-t)^{m/(1-m)}}\bigg|_{t=T(1-e^{-s})}\,,\qquad W_\infty(x)=v_{1/(1-m)}[z_\infty,\l_\infty](x)\,,
\]
so that \eqref{eq:vanishing} and \eqref{delpinosaez 1} imply, setting for short $w(s)=w(s,\cdot)$,
\begin{equation}
  \label{l2star convergence}
  \lim_{s\to +\infty}\|w(s)-W_\infty\|_{L^{2^\star}(\R^n)}=0\,.
\end{equation}
Recalling the notation $p=2^\star-1$, we see that $w$ satisfies
\begin{equation} \label{d_su_ds_wp}
\frac{d}{ds}w^p=\Delta w+\frac{1}{1-m}\,w^p\qquad \text{on }(0,\infty)\times \R^n\,,
\end{equation}
while of course $W=W_\infty$ is a solution to
\begin{equation}
\label{eq:FDbubble}
\Delta W+\frac{1}{1-m}\,W^p=0\qquad\text{on }\R^n\,.
\end{equation}
Let us consider the functional
$$
J[v]= \int_{\R^n}\frac{|\nabla v|^2}2 - \frac{1}{1-m}\int_{\R^n}\frac{v^{2^\star}}{2^\star}\,.
$$
By \eqref{d_su_ds_wp} we compute
\begin{equation}
\label{eq:diff J}
\frac{d}{ds}J[w(s)]=-\int_{\R^n}\Bigl(\Delta w(s)+\frac{1}{1-m}w(s)^p\Bigr)\frac{d}{ds}w(s)
=-\frac{1}{p}\int_{\R^n}\Bigl(\frac{\Delta w(s)}{w(s)^p}+\frac{1}{1-m}\Bigr)^2\,w(s)^{2^\star}\,,
\end{equation}
so that $s\mapsto J[w(s)]$ is decreasing, with $J[w(s)]\ge J[W_\infty]$ by Fatou's lemma and \eqref{l2star convergence}.
Exploiting the fact that
  $$
  \R\ni c\mapsto \int_{\R^n}\Bigl(\frac{\Delta w(s)}{w(s)^p}-c\Bigr)^2\,w(s)^{2^\star}
  $$
  attains its minimum at
  \[
  c=\frac{\int_{\R^n}\frac{\Delta w(s)}{w(s)^p}\,w(s)^{2^\star}}{\int_{\R^n}w(s)^{2^\star}}=-\frac{\int_{\R^n}|\nabla w(s)|^2}{\int_{\R^n}w(s)^{2^\star}}=-K_0(w(s))\,,
  \]
by H\"older inequality we find that
 \begin{eqnarray}
  \label{eq:bound 1delta}
  \delta(w(s))^2&\leq&
  \biggl(\int_{\R^n}w(s)^{2^\star}\biggr)^{2/n}\,\int_{\R^n}\Bigl(\frac{\Delta w(s)}{w(s)^p}+K_0(w(s))\Bigr)^2\,w(s)^{2^\star}
  \\\nonumber
  &\leq&
  C(n)\,\int_{\R^n}\Bigl(\frac{\Delta w(s)}{w(s)^p}+\frac1{1-m}\Bigr)^2\,w(s)^{2^\star}.
  \end{eqnarray}
Hence, thanks to \eqref{l2star convergence}, there exists $s_0$ (depending on the initial datum $u_0$) such that
\begin{equation}
  \label{fixed L2star norm}
  c(n)\le \int_{\R^n}w(s)^{2^\star}\le C(n)\qquad\forall\,s\ge s_0\,.
\end{equation}
Combining \eqref{eq:diff J}, \eqref{eq:bound 1delta}, and \eqref{fixed L2star norm} we find
\begin{equation}
  \label{de J}
  \de(w(s))^2\le -C(n)\,\frac{d}{ds}J[w(s)]\qquad\forall\,s\ge s_0\,,
\end{equation}
so that
\[
\int_{s_0}^\infty\de(w(s))^2\,ds\le C(n)\Big(J[w(s_0)]-J[w(\infty)]\Big)\le C(n)\Big(J[w(s_0)]-J[W_\infty]\Big)<\infty\,.
\]
In particular we can find a sequence $s_j\to\infty$ such that $\de(w(s_j))\to 0$ as $j\to\infty$. We can thus apply Struwe's theorem to $\{w(s_j)\}_{j\in\N}$. Since \eqref{l2star convergence} excludes bubbling, we conclude that
\[
\lim_{s\to\infty}J[w(s)]=\lim_{j\to\infty}J[w(s_j)]=J[W_\infty]\,.
\]
Furthermore, using again \eqref{l2star convergence}, this implies
\begin{equation}
  \label{limit dir}
  \lim_{s\to\infty}\int_{\R^n}|\nabla w(s)|^2=\int_{\R^n}|\nabla W_\infty|^2\,,
\end{equation}
and thus also
\begin{equation}
  \label{strong dir}
  \lim_{s\to\infty}\|\nabla w(s)-\nabla W_\infty\|_{L^2(\R^n)}=0\,,
  \qquad
  \lim_{s\to\infty}K_0(w(s))=\frac{\int_{\R^n}|\nabla W_\infty|^2}{\int_{\R^n}W_\infty^{2^\star}}=\frac1{1-m}\,.
\end{equation}
Now, for any $s>0$, denote by $W(s)$ the unique minimizer of
\begin{equation}
\label{eq:Ws min}
W\mapsto \|\nabla w(s)-\nabla W\|_{L^2(\R^n)}
\end{equation}
among all positive functions $W$ satisfying \eqref{eq:FDbubble}. Then, setting $w(s)=W(s)+\rho(s)$,
\begin{equation}
  \label{Ws and rhos}
  \qquad \int_{\R^n}\nabla W(s)\cdot\nabla\rho(s)=0=\int_{\R^n}W(s)^{p}\rho(s)\,,
  \qquad\lim_{s\to\infty}\|\nabla\rho(s)\|_{L^2(\R^n)}=0\,,
\end{equation}
where the last limit follow by \eqref{strong dir} and by the minimality property of $W(s)$.  Since $J$ is constant on solutions of \eqref{eq:FDbubble} (so $J[W(s)]=J[W_\infty]$), we can expand $J$ around $W(s)$ with the aid of \eqref{Ws and rhos}, to find out that
\begin{eqnarray}\nonumber
I[w(s)]&:=&J[w(s)]-J[W_\infty]=J[w(s)]-J[W(s)]
\\
&=&
\frac{1}{2} \int_{\R^n}|\nabla \rho(s)|^2 -
\frac{1}{(1-m)\,2^\star}\int_{\R^n}\Bigl(\bigl(W(s)+\rho(s)\bigr)^{2^\star}-W(s)^{2^\star} - 2^\star \,W(s)^p\,\rho(s)\Bigr)\nonumber \\
&\approx& \int_{\R^n}|\nabla \rho(s)|^2 - \frac{p}{1-m}\int_{\R^n}W(s)^{p-1}\rho(s)^2\label{eq:I expansion}\\
&\approx& \int_{\R^n}|\nabla \rho(s)|^2=\int_{\R^n}|\nabla w(s)-\nabla W(s)|^2\qquad \forall\, s \geq s_0\,,\nonumber
\end{eqnarray}
where  $a\approx b$ means that $a/C(n) \leq b \leq C(n) \,a$ for some positive dimensional constant $C(n)$,
 the second $\approx$ follows from the spectral gap estimate \eqref{spectral gap lowerbound}.

We now notice that by \eqref{limit dir}, up to possibly increase the value of $s_0$, we can apply Theorem \ref{thm main} to each $w(s)$ with $s\ge s_0$: in particular, for every $s\ge s_0$ there exists $\overline{W}(s)$ such that
\begin{equation}
  \label{our}
-\Delta \overline{W}(s)=K_0(w(s))\,\overline{W}(s)^p\,,
\qquad \|\nabla w(s)-\nabla \overline{W}(s)\|_{L^2(\R^n)}\leq C(n)\,\delta(w(s))\,.
\end{equation}
By \eqref{fixed L2star norm} and by \eqref{eq:diff J},
  \begin{eqnarray}\nonumber
  \Bigl|K_0(w(s))-\frac{1}{1-m}\Bigr|^2&=&\bigg|\frac{\int_{\R^n}\bigl(\frac{\Delta w(s)}{w(s)^p}+\frac{1}{1-m}\bigr)\,w(s)^{2^\star}}{\int_{\R^n}w(s)^{2^\star}}\bigg|^2
  \\\label{eq:bound 2}
  &\leq& C(n)\, \int_{\R^n}\Big(\frac{\Delta w(s)}{w(s)^p}+\frac{1}{1-m}\Big)^{2}w(s)^{2^\star}=-C(n)\,\frac{d}{ds}\,J[w(s)]\,.\hspace{1cm}
    \end{eqnarray}
Setting $\alpha(s)=[(1-m)\,K_0(w(s))]^{1/(p-1)}$, the function $\hat W(s)=\alpha(s)\,\overline W(s)$ satisfies \eqref{eq:FDbubble}, with
$$
\|\nabla \hat W(s)-\nabla \overline{W}(s)\|_{L^2(\R^n)}\leq C(n)\,\Bigl|K_0(w(s))-\frac{1}{1-m}\Bigr|\,,
$$
and hence, by triangle inequality, \eqref{our}, \eqref{de J}, and \eqref{eq:bound 2},
$$
\|\nabla w(s)-\nabla \hat{W}(s)\|_{L^2(\R^n)}^2
\leq C(n)\,\Bigl(\delta(w(s))^2+\Bigl|K_0(w(s))-\frac{1}{1-m}\Bigr|^2\Bigr)\le -C(n)\,\frac{d}{ds}J[w(s)]\,.
$$
By the minimality property of $W(s)$ and since $J[w(s)]$ and $I[w(s)]$ differ by a constant, we conclude that
\[
\|\nabla\rho(s)\|_{L^2(\R^n)}^2=\|\nabla w(s)-\nabla W(s)\|_{L^2(\R^n)}^2\le -C(n)\,\frac{d}{ds}I[w(s)]
\]
for every $s\ge s_0$, and thus, by \eqref{eq:I expansion},
$$
\frac{d}{ds}I[w(s)]\leq -C(n)\,I[w(s)]\qquad \forall\, s \geq s_0\,.
$$
This proves that, for a constant $C_*$ depending on the initial datum $u_0$,
\begin{equation}
  \label{iii}
  I[w(s)] \leq C_*\,e^{-\k(n)\,s}\qquad \forall\,s\ge0
\end{equation}
(note that the above inequality is trivially true on $[0,s_0]$
by choosing $C_*$ large enough)
and hence, applying \eqref{eq:I expansion} again, we deduce that
\begin{equation} \label{trump}
\int_{\R^n}|\nabla w(s)-\nabla W(s)|^2 \leq  C_*\,e^{-\k(n)\,s}\qquad \forall\, s\geq 0\,.
\end{equation}
In other words, we have prove that for all $s$ there exists a function $W(s)$ solving \eqref{eq:FDbubble} which is exponentially close to $w(s)$. To conclude the result, we need to show that $W(s)$ is exponentially close to $W_\infty$.
To this aim
we notice that, from \eqref{d_su_ds_wp},
$$
\frac{d}{ds} w^{2^\star}=\frac{p+1}{p} w  \frac{d}{ds} w^{p} = \frac{2n}{n+2} \left(\Delta w + \frac{1}{1-m} w^p \right) w,
$$
so we obtain
$$
|w(s)-w(t)|^{2^\star} \leq |w(s)^{2^*} - w(t)^{2^*}| \leq C(n) \Big{|} \int_t^s \left(\frac{\Delta w(r)}{w(r)^p} + \frac{1}{1-m} \right) w(r)^{2^*}  dr \Big{|} \,,
$$
and H\"older inequality yields
\begin{equation} \label{cauchy}
\int_{\R^n} |w(s)-w(t)|^{2^*} \leq C(n) \int_t^s \delta(w(r)) dr \ \quad \forall\,s>t>0 \,.
\end{equation}
Now, by \eqref{de J} and \eqref{iii}
\[
\Big(\int_k^{k+1}\de(w(r))\,dr\Big)^2\le\int_k^{k+1}\de(w(r))^2\le C(n)\,\big(I[w(k)]-I[w(k+1)]\big)\le C_*\,e^{-\k(n)\,k}
\]
which implies
$$
\int_t^\infty \delta(w(s))\,ds \leq C_* e^{-c(n) t} \quad \quad \forall\,t >0 \,.
$$
Hence from \eqref{cauchy} we obtain
$$
\int_{\R^n} |w(s)-w(t)|^{2^*} \leq C_* e^{-c(n) t} \quad \quad \forall\,s > t >0 \,.
$$
From Sobolev inequality and \eqref{trump}, it follows that $\{W(s)\}_{s>0}$ is a Cauchy family in $L^{2^\star}(\R^n)$ satisfying the exponential bound
$$
\int_{\R^n} |W(s)-W(t)|^{2^*} \leq C_* e^{-c(n) t} \quad \quad \forall\,s > t >0 \,.
$$
Since $\|\nabla W(s)-\nabla W_\infty\|_{L^2(\R^n)}\to 0$
as $s \to \infty$ (see \eqref{strong dir} and \eqref{trump}),
we conclude that 
$$
\int_{\R^n} |W_\infty -W(t)|^{2^*} \leq C_* e^{-c(n) t} \quad \quad \forall\,t >0 \,,
$$
that combined with \eqref{trump} and Sobolev inequality yields
\begin{equation}
  \label{eq:exp w}
  \int_{\R^n} |W_\infty -w(t)|^{2^*} \leq C_* e^{-c(n) t} \quad \quad \forall\,t >0 \,.
\end{equation}
Now, to conclude the proof, we argue as follows:
let $F:\mathbb S^n\to \R^n$ denote the
inverse stereographic projection (see \eqref{eq:F}),
and let $v(s):\mathbb{S}^n\to \R$ be defined from $w(s)$ via the transformation
\begin{equation}
  \label{transformation}
  w(s,x)=\Bigl(\frac{2}{1+|x|^2}\Bigr)^{(n-2)/2}v(s,F(x))\,.
\end{equation}
Then $v(s)$ solves the equation
\begin{equation}
\label{eq:v p}
\frac{d}{ds}v^p=\Delta_{\mathbb{S}^n}v+\frac{1}{1-m}v^p -\frac{n(n-2)}{4}\,v\qquad\mbox{on $(0,\infty)\times\mathbb{S}^n$}
\end{equation}
(see for instance \cite{YYLiPARTI} and \cite[Equation (2.3)]{DPS}) and \eqref{eq:exp w} translates into
\begin{equation}
\label{eq:exp v}
\int_{\mathbb{S}^n}| v(s)-v_\infty|^{2^\star} \leq  C_*\,e^{-\k(n)\,s}\qquad \forall\, s\geq 0\,,
\end{equation}
where $v_\infty$ is the stationary solution of \eqref{eq:v p} corresponding to $W_\infty$ under the transformation \eqref{transformation}.
Since $v(s)$ is uniformly bounded away from zero and infinity for $s$ large (see \cite[Proposition 5.1]{DPS}), it follows by  \eqref{eq:exp v} and parabolic regularity that (up to replacing $\k(n)$ by $\k(n)/2$)
$$
\|v(s)- v_\infty\|_{L^\infty(\mathbb{S}^n)} \leq  C_*\,e^{-\k(n)\,s}\qquad \forall\, s\geq 0\,.
$$
Going back to the original variables,
this implies that \eqref{eq:to prove FD} holds, concluding the proof.

\bibliography{references}
\bibliographystyle{is-alpha}

\end{document}